\documentclass{article}
\usepackage[utf8]{inputenc}
\usepackage[cm]{fullpage}

\usepackage{booktabs}
\usepackage{amsmath, amssymb, amsthm}
\usepackage{multicol}
\usepackage{multirow}
\usepackage{tablefootnote}
\usepackage{xspace}
\usepackage[shortlabels, inline]{enumitem}
\usepackage{thmtools}

\usepackage{hyperref}
\hypersetup{
  colorlinks=true,
  citecolor=black,
  filecolor=black,
  linkcolor=blue,
  urlcolor=red
}
\usepackage{zref-clever}
\zcsetup{nameinlink=true, capfirst}
\newcommand{\cref}[1]{\zcref[S]{#1}}

\setenumerate{noitemsep, topsep=0pt}

\theoremstyle{plain}
\declaretheorem[name=Theorem, numberwithin=section]{theorem}
\declaretheorem[name=Corollary, sibling=theorem]{cor}

\theoremstyle{definition}
\declaretheorem[name=Definition, sibling=theorem]{defn}

\zcRefTypeSetup{theorem}{
  Name-sg = Theorem,
  name-sg = Theorem,
  Name-pl = Theorems,
  name-pl = theorems,
}

\zcRefTypeSetup{cor}{
  Name-sg = Corollary,
  name-sg = Corollary,
  Name-pl = Corollaries,
  name-pl = corollaries,
}

\zcRefTypeSetup{defn}{
  Name-sg = Definition,
  name-sg = Definition,
  Name-pl = Definitions,
  name-pl = definitions,
}

\usepackage[bibencoding=utf8, giveninits=true, sortcites]{biblatex}
\newcommand{\GAP}{\textsc{GAP}~\cite{GAP4}\xspace}
\newcommand{\Smallsemi}{\textsc{Smallsemi}~\cite{Smallsemi}\xspace}
\newcommand{\Semigroups}{\textsc{Semigroups}~\cite{Semigroups}\xspace}
\newcommand{\Semirings}{\textsc{Semirings}~\cite{Semirings}\xspace}
\newcommand{\bliss}{\textsc{bliss}~\cite{bliss, junttila2007}\xspace}
\newcommand{\alg}{\textsc{alg}~\cite{baueralg}\xspace}
\newcommand{\Sym}{\operatorname{Sym}}
\newcommand{\Aut}{\operatorname{Aut}}

\title{Counting finite semirings}
\author{J. Edwards, J. D. Mitchell, and P. Ragavan}
\date{\today}

\begin{document}

\maketitle

\begin{abstract}
  In this short note we count the finite semirings up to isomorphism,
  and up to isomorphism or anti-isomorphism for some small values of
  $n$; for which we utilise the existing library of small semigroups
  in the \GAP package \Smallsemi.
\end{abstract}

\section{Introduction}

Enumeration of algebraic and combinatorial structures of finite order up to
isomorphism is a classical topic. Among the algebraic structures considered are
groups~\cite{BESCHE2002,A000001}, rings~\cite{Blackburn2022, Fine1993,
Kruse1970, A027623}, near-rings~\cite{Chow2024, SONATA, A305858},
semigroups and monoids~\cite{Distler2010, Distler2012, Distler2013,
  Forsythe1955, Grillet1996, Grillet2014, Jrgensen1977, Klee1958, Motzkin1955,
Plemmons1967, A027851, A058129, Satoh1994, A001426, A001423},
inverse semigroups and monoids~\cite{Malandro2019, A001428, A234843,
A234844, A234845}, and many more too numerous to mention\footnote{The disparity
  in the number of references for semigroups and monoids and the other
  algebraic structures is a consequence of the authors' familiarity with the
  literature for semigroups and monoids, and there are likely many other
  references that could have been included were it not for us not knowing about
them.}. In this short note we count the number of finite semirings up to
isomorphism and up to isomorphism or anti-isomorphism for $n \leq 6$; see
\cref{tab:semirings,tab:comm-semirings,tab:ai-semirings,tab:comm-ai-semirings}.
We also count several special classes of semirings for (slightly)
larger values of $n$.

This short note was initiated by an email from M. Volkov to the second author
in February of 2025 asking if it was possible to verify the claim
in~\cite{Ren2025} that the number of ai-semirings (additively
idempotent) up to isomorphism with $4$
elements is $866$ (see the caption of \cref{tab:ai-semirings} for the
definition). After some initial missteps it was relatively
straightforward to verify that this number is correct, by using the library of
small semigroups in the \GAP package \Smallsemi. This short note arose from
these first steps. In contrast to groups or rings, where the numbers of
non-isomorphic objects of order $n$ are known for relatively large values of
$n$, the number of semigroups of order $10$ (up to isomorphism) is apparently
not known exactly (although the number up to isomorphism or anti-isomorphism
is known~\cite{Distler2012}). The paper \cite{Kleitman1976} purports to show
that almost all semigroups of order $n$ (as $n$ tends to infinity) up to
isomorphism or anti-isomorphism are $3$-nilpotent (and indeed $99.4\%$ of the
semigroups of order $8$ are $3$-nilpotent), and it is shown
in~\cite{Distler2012ab} that the number of $3$-nilpotent semigroups of order
$10$ is approximately $10 ^ {19}$. This number is likely close to the exact
value. Perhaps unsurprisingly, from the perspective of counting up to
isomorphism, it seems that semirings have more in common with semigroups than
with rings or groups. Very roughly speaking, rings and groups are highly
structured, providing strong constraints that facilitate their enumeration. On
the other hand, semigroups, and seemingly semirings also, are less structured,
more numerous, and harder to enumerate.

We begin with the definition of a semiring, which is a natural generalisation
of the notion of a ring.

\begin{defn}[Semiring]
  \label{def:semiring}
  A \emph{semiring} is a set \(S\) equipped with two binary
  operations \(+\) and  \(\times\) such that:
  \begin{enumerate}
    \item \((S, +)\) is a commutative semigroup ($(x + y) + z = x + (y + z)$
      and $x + y = y + x$ for all $x, y, z\in S$);
    \item \((S, \times)\) is a semigroup ($x\times (y \times z) = (x\times
      y)\times z$ for all $x,y, z\in S$); and
    \item multiplication $\times$ distributes over addition $+$
      ($x \times (y + z) = (x \times y) + (x \times z)$ and
      $(y + z) \times x = (y \times x) + (z \times x)$ for all $x, y, z \in S$).
  \end{enumerate}
\end{defn}

We note that some authors require that $(S, +)$ is a commutative monoid with
(additive) identity $0$ such that $0 \times x = x \times 0 = 0$ for
all $x\in S$; see, for example, \cite{Lothaire2005, Sakarovitch2009}. We
do not add this requirement, and refer to such objects as \textit{semirings
with zero}. The numbers of semirings with zero are discussed
in~\cite{stackexchange} and in~\cite{jipsen}; see also~\alg which we will
discuss a little further below.
An \textit{ai-semiring} is a semiring $S$ where $x + x = x$ for all $x \in S$.
In~\cite{Zhao2020} it was shown that there are
$61$ ai-semirings of order $3$.

Recall that if $(S, \times)$ and $(T, \otimes)$ are semigroups, then $\phi:S
\to T$ is a \textit{(semigroup) homomorphism} if $(x\times y)\phi =
(x)\phi\otimes (y)\phi$ for all $x, y\in S$. Note that we write mappings to the
right of their arguments and compose them from left to right. If $\phi: S \to
T$ is a semigroup homomorphism and $\phi$ is bijective, then $\phi$ is an
\textit{(semigroup) isomorphism}. A semigroup isomorphism from a semigroup
$(S, \times)$ to itself is called an \textit{automorphism}
and the group of all such automorphisms is denoted $\Aut(S, \times)$.

In this short note we are concerned with counting semirings up to isomorphism,
and so our next definition is that of an isomorphism.

\begin{defn}[Semiring isomorphism]\label{defn-semiring-iso}
  We say that two semirings \((S,+,\times)\) and \((T, \oplus, \otimes)\) are
  \emph{isomorphic} if there exists a bijection \(\phi: S \to T\) which is
  simultaneously a semigroup isomorphism from $(S, +)$ to $(T, \oplus)$ and
  from $(S, \times)$ to $(T, \otimes)$. We refer to $\phi$ as a
  \textit{(semiring) isomorphism}.
\end{defn}

Since a semiring is comprised of two semigroups, enumerating semirings is
equivalent to enumerating those pairs consisting of an additive semigroup $(S,
+)$ and a multiplicative semigroup $(S, \times)$ such that $\times$ distributes
over $+$. The next theorem indicates which $(S, \times)$ we should consider
for each of the additive semigroups $(S, +)$.

We denote the symmetric group on the set $S$ by $\Sym(S)$. If $\sigma\in
\Sym(S)$ and $\cdot: S \times S \to S$ is a binary operation, then we define
the binary operation $\cdot ^ \sigma: S\times S \to S$ by
\begin{equation}\label{equation-action}
  x \cdot ^ \sigma y = ((x)\sigma^{-1} \cdot (y)\sigma ^ {-1})\sigma.
\end{equation}
It is straightforward to verify that \eqref{equation-action} is a (right, group)
action of $\Sym(S)$ on the set of all binary operations on $S$. Clearly if
$\cdot$ is associative, then so too is $\cdot ^ \sigma$ for every $\sigma \in
\Sym(S)$. The group of automorphisms $\Aut(S,\cdot)$ of a
semigroup $(S, \cdot)$ coincides with the stabiliser of the operation $\cdot$
under the action of $\Sym(S)$ defined in \eqref{equation-action}.

Recall that if $H$ and $K$ are subgroups of a group $G$, then the
\textit{double cosets} $H\backslash G / K$ are the sets of the form
$\{hgk\mid h\in H, k \in K\}$ for $g\in G$.
The next theorem is key to our approach for counting semirings.

\begin{theorem}
  \label{thm:isomorphism-condition}
  Let \((S, +)\) be a commutative semigroup, let \((S, \times)\) be a
  semigroup, and let \(\sigma, \tau\in\Sym(S)\). Then there exists a
  semigroup isomorphism $\phi\in \Aut(S, +)$
  from \((S,\times^\sigma)\) to \((S, \times ^ \tau)\) if
  and only if \(\sigma\) and \(\tau\) belong to the same double coset in
  \(\Aut(S, \times) \backslash \Sym(S) / \Aut(S, +)\).
\end{theorem}

\begin{proof}
  ($\Rightarrow$) Suppose that $\phi\in \Aut(S, +)$ is a semigroup isomorphism
  from $(S, \times ^ \sigma)$ to $(S, \times^\tau)$. Then $(x \times ^ \sigma
  y)\phi = (x)\phi\times ^ \tau (y)\phi$ for all $x, y\in S$. If $x, y\in S$
  are arbitrary, then
  \[
    ((x)\sigma ^ {-1} \times (y)\sigma ^ {-1})\sigma \phi \tau ^ {-1} =
    (x \times ^ \sigma y)\phi\tau ^ {-1}
    = ((x)\phi \times ^ {\tau} (y)\phi)\tau ^ {-1}
    = ((x)\phi\tau^{-1}\times (y)\phi\tau^{-1})\tau\tau^{-1}
    = (x)\phi\tau ^ {-1} \times (y)\phi\tau ^ {-1}.
  \]
  If we set $\gamma = \sigma\phi\tau^{-1}$, $p =
  (x)\sigma^{-1}$, and $q=(y)\sigma^{-1}$, then
  $(p\times q)\gamma = (p)\gamma \times (q)\gamma$
  and so $\gamma, \gamma ^ {-1} \in \Aut(S, \times)$.
  Rearranging we obtain $\tau = \gamma^{-1}\sigma\phi$. Since $\phi \in
  \Aut(S, +)$, we conclude that $\tau$ and $\sigma$ lie in the
  same double coset of $\Aut(S, \times) \backslash \Sym(S) / \Aut(S, +)$.

  ($\Leftarrow$) Suppose that $\sigma$ and $\tau$ are in the
  same double coset of $\Aut(S, \times) \backslash \Sym(S) / \Aut(S, +)$.
  Then there exists $\alpha\in\Aut(S, \times)$ and $\beta\in\Aut(S, +)$ such
  that $\tau=\alpha\sigma\beta$.

  We will show that $\beta$ is a semigroup isomorphism from $(S, \times
  ^{\sigma})$ to $(S, \times^{\tau})$:
  \begin{align*}
    (x\times^\sigma y)\beta
    &= (x\sigma ^{-1} \times y\sigma ^{-1})\sigma \beta
    = (x\sigma ^{-1} \times y\sigma ^{-1})\alpha^{-1}\tau
    = (x\sigma ^{-1}\alpha^{-1} \times y\sigma ^{-1}\alpha^{-1})\tau &&
    \alpha^{-1} \in \Aut(S, \times)\\
    &= (x\beta\tau^{-1} \times y\beta\tau^{-1})\tau
    = (x\beta \times^{\tau} y\beta).\qedhere
  \end{align*}
\end{proof}

An immediate corollary of \cref{thm:isomorphism-condition} is the following.

\begin{cor}\label{cor-double-coset-reps-iso}
  Let \((S, +)\) be a commutative semigroup, let \((S, \times)\) be a
  semigroup, and let \(\sigma\) and \(\tau\) belong to the
  same double coset in
  \(\Aut(S, \times) \backslash \Sym(S) / \Aut(S, +)\).
  Then the following hold:
  \begin{enumerate}
    \item
      \((S, +, \times^\sigma)\) is a semiring if and only if
      \((S, +, \times^\tau)\) is a semiring;
    \item if \((S, +, \times^\sigma)\) and \((S, +, \times^\tau)\) are
      semirings, then they are isomorphic.
  \end{enumerate}
\end{cor}

If $(S, \times)$ and $(S, \otimes)$ are semigroups, then $(S,
\times)$ is said
to be \textit{anti-isomorphic} to $(S, \otimes)$ if there exists a bijection
$\phi: S \to S$ such that $(x\times y)\phi = (y)\phi\otimes (x)\phi$ for all
$x, y\in S$.
The bijection $\phi$ is referred to as an \textit{anti-isomorphism}. If the
operations $\times$ and $\otimes$ coincide, then $\phi$ is an
\textit{anti-automorphism}. Clearly the composition of two
anti-automorphisms
is an automorphism, and the composition of an anti-automorphism and an
automorphism is an anti-automorphism. As such the set of all
automorphisms or
anti-automorphisms forms a group under composition of functions; we
denote this
group by $\Aut^*(S, \times)$.

Similarly, using the obvious analogue of \cref{defn-semiring-iso}, we can
define anti-isomorphic semirings. It is routine to adapt the proof of
\cref{thm:isomorphism-condition} to prove the following.

\begin{theorem}
  \label{cor:equiv-condition}
  Let \((S, +)\) be a commutative semigroup, let \((S, \times)\) be a
  semigroup, and let \(\sigma, \tau\in\Sym(S)\).
  Then there exists a semigroup isomorphism or
  anti-isomorphism \(\phi\in \Aut(S, +)\)  from \((S, \times ^
  \sigma)\) to \((S, \times^\tau)\)
  if and only if \(\sigma\) and \(\tau\) belong to the same double coset in
  \(\Aut^*(S, \times) \backslash \Sym(S) / \Aut(S, +)\).
\end{theorem}

We also obtain a corollary of \cref{cor:equiv-condition} analogous to
\cref{cor-double-coset-reps-iso}.

\begin{cor}\label{cor-double-coset-reps-equiv}
  Let \((S, +)\) be a commutative semigroup, let \((S, \times)\) be a
  semigroup, and let \(\sigma\) and \(\tau\) belong to the
  same double coset in
  \(\Aut^*(S, \times) \backslash \Sym(S) / \Aut(S, +)\).
  Then the following hold:
  \begin{enumerate}
    \item
      \((S, +, \times^\sigma)\) is a semiring if and only if
      \((S, +, \times^\tau)\) is a semiring;
    \item if \((S, +, \times^\sigma)\) and \((S, +, \times^\tau)\) are
      semirings, then they are isomorphic or anti-isomorphic.
  \end{enumerate}
\end{cor}

\section{Tables of results}

The approach we implemented in the \GAP package \Semirings to compute the
numbers in this section is as follows. It follows from
\cref{cor-double-coset-reps-iso,cor-double-coset-reps-equiv} that to
enumerate the semirings on a set $S$ up to isomorphism (or up to
isomorphism or anti-isomorphism) it suffices to consider every
commutative semigroup $(S, +)$ and semigroup $(S, \times)$ and to check
distributivity of $+$ and $\times^{\sigma}$ for one representative
$\sigma$ of every double coset in \(\Aut(S, \times) \backslash \Sym(S) /
\Aut(S, +)\) (or in \(\Aut^*(S, \times) \backslash \Sym(S) / \Aut(S,
+)\)). The semigroups with order at most $8$ up to isomorphism and
anti-isomorphism are available in \Smallsemi.  Since $S$ is small, it is
relatively straightforward to compute $\Aut(S, \times)$ and $\Aut(S, +)$:
the problem of computing $\Aut(S, \times)$ reduces to that of
computing the automorphism group of a graph associated to $(S, \times)$;
see \cite{Miller1977} for details. The automorphisms of $(S, +)$ can also
be computed in this way. There are several software tools for computing
automorphism groups of graphs, for example, \bliss. This is the approach
used by the \Semigroups package for \GAP, which the authors used to
compute $\Aut(S, \times)$ and $\Aut(S, +)$ in this paper. Finally, \GAP
contains functionality for computing double coset representatives; see
\cite[Section 4.6.8]{Holt2005} for details.

Below are some tables listing the numbers of semirings
(\cref{def:semiring}) up
to isomorphism, and up to isomorphism or anti-isomorphism, with various
properties for some small values of $n\in \mathbb{N}$. As far as
we know, many
of the numbers in these tables were not previously known. In
particular, we are
not aware of any results in the literature about the number of
semirings up to
isomorphism or anti-isomorphism. Some of the numbers in the
tables below can
be found using \alg, although using \alg is considerably slower
than the approach described here, largely because the precomputed data for
small semigroups available in \Smallsemi does a lot of the heavy lifting.

As a sanity check, where possible, we have checked the numbers produced by
\Semirings against those in the literature, and those we could produce using
\alg. With one exception, these numbers agreed with each other.
\alg computes
that there are $57,443$ semirings (using \cref{def:semiring}) of
order $5$ up
to isomorphism, whereas \Semirings computes $57,427$. Using \GAP
and \Semirings
we attempted to verify whether or not the output of \alg was correct.
Unfortunately, it appears that $16$ of the pairs of multiplication tables
output by \alg do not satisfy one or the other of the
distributivity conditions
from \cref{def:semiring}(c); see \cite{bauer_alg_issue16} for more details.
This precisely accounts for the difference between the number output by
\Semirings and the number output by \alg.

The numbers present in the tables below were computed using a variety of
desktop and laptop computers with the following specs: a 2025 MacBook Pro M4
Pro with 14 threads and 48GB of RAM; an Intel(R) Core(TM) i7-8700
CPU @ 3.20GHz
processor with 12 threads and 64GB of RAM; an Intel(R) Core(TM)
i7-12700KF @ 5GHz with 20 threads and 64GB of RAM. We mostly include this
information to indicate that we did not make use of particularly powerful
computers. The longest of the computation in the tables below
took approximately $2400$ CPU hours. This was the time taken to count the
semirings with $0$ up to isomorphism or anti-isomorphism with $7$ elements.

These tables are only a sample of the results that can be obtained
using \Semirings.

\begin{table}[ht]
  \centering
  \begin{tabular}{l|r|r|r|r|r|r|r|r}
    $n$
    & \multicolumn{4}{c|}{up to isomorphism}
    & \multicolumn{4}{c}{up to isomorphism or anti-isomorphism} \\
    \midrule
    & \multicolumn{1}{c|}{\shortstack{no additional\\constraints}} &
    \multicolumn{1}{c|}{with $0$}
    & \multicolumn{1}{c|}{with $1$} & \multicolumn{1}{c|}{with $0$ + $1$}
    & \multicolumn{1}{c|}{\shortstack{no additional\\constraints}} &
    \multicolumn{1}{c|}{with $0$}
    & \multicolumn{1}{c|}{with $1$} & \multicolumn{1}{c}{with $0$ + $1$}
    \\
    \midrule
    1   & 1         & 1         & 1          & 1      & 1               & 1         & 1         & 1      \\
    2   & 10        & 4         & 4          & 2      & 9               & 4         & 4         & 2      \\
    3   & 132       & 22        & 22         & 6      & 106             & 20        & 21        & 6      \\
    4   & 2,341     & 283       & 169        & 40     & 1,713           & 226       & 155       & 38     \\
    5   & 57,427    & 4,717     & 1,819      & 295    & 38,247          & 3,365     & 1,561     & 262    \\
    6   & 7,571,579 & 108,992   & 41,104     & 3,246  & 4,102,358       & 71,138    & 30,112    & 2,681  \\
    7   & -         & 8,925,672 & 11,679,328 & 59,314 & 48,152,448,707  & 4,910,824 & 6,268,858 & 43,331 \\
  \end{tabular}
  \caption{Numbers of semirings (\cref{def:semiring}) with $n$
    elements up to
    isomorphism and up to isomorphism or anti-isomorphism. See
    \cite{MSsemirings,MSsemiringsWithOneAndZero,MSsemiringsWithOne}
    for some of
  these numbers up to isomorphism.}\label{tab:semirings}
\end{table}

\begin{table}[ht]
  \centering
  \begin{tabular}{l|r|r|r|r}
    $n$
    & \multicolumn{4}{c}{up to isomorphism}\\
    \midrule
    & \multicolumn{1}{c|}{\shortstack{no additional\\constraints}} &
    \multicolumn{1}{c|}{with $0$}
    & \multicolumn{1}{c|}{with $1$} & \multicolumn{1}{c}{with $0$ + $1$}
    \\
    \midrule
    1   & 1           & 1       & 1       & 1       \\
    2   & 8           & 4       & 4       & 2       \\
    3   & 80          & 18      & 20      & 6       \\
    4   & 1,067       & 169     & 141     & 36      \\
    5   & 18,188      & 1,990   & 1,276   & 228     \\
    6   & 543,458     & 32,212  & 17,621  & 2,075   \\
    7   & 162,744,745 & 799,354 & 690,924 & 25,640  \\
    8   & -           & -       & -       & 791,061 \\
  \end{tabular}
  \caption{Numbers of commutative semirings (i.e. those satisfying
  $x\times y = y \times x$ for all $x,y\in S$) with $n$ elements.}
  \label{tab:comm-semirings}
\end{table}

\begin{table}[ht]
  \centering
  \begin{tabular}{l|r|r|r|r|r|r|r|r}
    $n$
    & \multicolumn{4}{c|}{up to isomorphism}
    & \multicolumn{4}{c}{up to isomorphism or anti-isomorphism} \\
    \midrule
    & \multicolumn{1}{c|}{\shortstack{no additional\\constraints}} &
    \multicolumn{1}{c|}{with $0$}
    & \multicolumn{1}{c|}{with $1$} & \multicolumn{1}{c|}{with $0$ + $1$}
    & \multicolumn{1}{c|}{\shortstack{no additional\\constraints}} &
    \multicolumn{1}{c|}{with $0$}
    & \multicolumn{1}{c|}{with $1$} & \multicolumn{1}{c}{with $0$ + $1$}
    \\
    \midrule
    1   & 1         & 1       & 1       & 1       & 1         & 1       & 1       & 1       \\
    2   & 6         & 2       & 2       & 1       & 5         & 2       & 2       & 1       \\
    3   & 61        & 12      & 11      & 3       & 45        & 10      & 10      & 3       \\
    4   & 866       & 129     & 73      & 20      & 581       & 93      & 64      & 18      \\
    5   & 15,751    & 1,852   & 703     & 149     & 9,750     & 1,207   & 574     & 125     \\
    6   & 354,409   & 33,391  & 9,195   & 1,488   & 205,744   & 20,142  & 6,835   & 1,150   \\
    7   & 9,908,909 & 729,629 & 164,163 & 18,554  & 5,470,437 & 415,527 & 109,880 & 13,171  \\
    8   & -         & -       & -       & 295,292 & -         & -       & -       & 116,274 \\
  \end{tabular}
  \caption{Numbers of ai-semirings (i.e. those satisfying $x + x
      = x$ for all
    $x\in S$, where ``ai'' stands for ``additively idempotent'') with $n$
    elements.
  }\label{tab:ai-semirings}
\end{table}

\begin{table}[ht]
  \centering
  \begin{tabular}{l|r|r|r|r}
    $n$
    & \multicolumn{4}{c}{up to isomorphism}\\
    \midrule
    & \multicolumn{1}{c|}{\shortstack{no additional\\constraints}} &
    \multicolumn{1}{c|}{with $0$}
    & \multicolumn{1}{c|}{with $1$} & \multicolumn{1}{c}{with $0$ + $1$}
    \\
    \midrule
      1   & 1          & 1         & 1       & 1      \\
      2   & 4          & 2         & 2       & 1      \\
      3   & 29         & 8         & 9       & 3      \\
      4   & 289        & 57        & 55      & 16     \\
      5   & 3,589      & 550       & 437     & 100    \\
      6   & 53,661     & 6,639     & 4,296   & 794    \\
      7   & 949,843    & 96,264    & 52,043  & 7,493  \\
      8   & 20,054,643 & 1,639,905 & 764,329 & 84,961 \\
  \end{tabular}
  \caption{Numbers of commutative ai-semirings (i.e. those satisfying
    $x\times y = y \times x$ and $x + x = x$ for all $x,y\in S$) with
  $n$ elements.}
  \label{tab:comm-ai-semirings}
\end{table}

\section*{Acknowledgements}
The authors would like to thank the anonymous referees for their helpful
comments.

\printbibliography

\end{document}